\documentclass[12pt]{article}

\usepackage[T1]{fontenc}
\usepackage[english]{babel}
\usepackage[a4paper,margin=1in]{geometry}
\usepackage{amsmath,amssymb,amsthm,mathtools}
\usepackage{microtype}
\usepackage{hyperref}
\usepackage{cleveref}
\allowdisplaybreaks
\hypersetup{
  hidelinks,
  pdftitle={Exact L2 Norms for Graph Intersections and Discrepancy of (n/2)-Uniform Hypergraphs},
  pdfauthor={Zhanping Yang},
  pdfsubject={Exact graph intersection second moments and half-dense uniform hypergraph discrepancy},
  pdfkeywords={graph intersection, discrepancy, L2 norm, degree variance, quasi-complete graph, quasi-star graph, uniform hypergraph}
}

\newtheorem{theorem}{Theorem}[section]
\newtheorem{lemma}[theorem]{Lemma}
\newtheorem{problem}[theorem]{Problem}

\theoremstyle{remark}
\newtheorem{remark}[theorem]{Remark}

\crefname{theorem}{theorem}{theorems}
\crefname{lemma}{lemma}{lemmas}
\crefname{corollary}{corollary}{corollaries}
\crefname{remark}{remark}{remarks}

\title{Graph Intersections and Hypergraph Discrepancy}
\author{
  Zhanping Yang\thanks{Corresponding author. Center for Discrete Mathematics, Fuzhou University, Fuzhou, Fujian 350108, China. Supported by the National Key R\&D Program of China (Grant No. 2023YFA1010202) and the National Natural Science Foundation of China (Grant No. 12371342). Email: \texttt{yangzp@fzu.edu.cn}.}
}
\date{}

\begin{document}
\maketitle

\begin{abstract}
Let $G$ and $H$ be graphs on $n$ vertices with $e(G)=p\binom{n}{2}$ and $e(H)=q\binom{n}{2}$,
where $p,q\in(0,1)$ are fixed. Bollob\'as and Scott studied intersections of two graphs under random relabelling and asked for the minimum $L^2$ norm of the deviation of their intersection size from its expected value. We answer this question by proving that, for graphs of order
$n$ with fixed edge densities $p,q\in(0,1)$, the minimum is
\begin{equation*}
\left(\sqrt{\frac{p(1-p)q(1-q)}{2}\cdot\frac{1-\sqrt R}{1+\sqrt R}}\right)n+o(n),
\end{equation*}
where $R=\max\{p,1-p,q,1-q\}$. The coefficient of the linear term is best possible.

Bollob\'as and Scott defined the positive and negative discrepancies of two hypergraphs as the largest deviations of their intersection size above and below its expected value under relabelling, respectively. They asked for the minimum product of the positive and negative discrepancies of two $(n/2)$-uniform hypergraphs on $n$ vertices with $\frac{1}{2}\binom{n}{n/2}$ edges. We answer this question by proving that the minimum product of these discrepancies is $0$ if $4\mid\binom{n}{n/2}$ and $1/4$ otherwise.
\end{abstract}



\section{Introduction}\label{sec:introduction}
The distribution of edges in graphs has long been studied through discrepancy and related extremal problems. Erd\H{o}s and Spencer~\cite{ErdosSpencer1971} proved that every graph $G$ on $n$ vertices contains a set $S\subseteq V(G)$ such that
\begin{equation*}
\left|e_G(S)-\frac12\binom{|S|}{2}\right|\ge cn^{3/2},
\end{equation*}
where $c>0$ is an absolute constant. Erd\H{o}s, Goldberg, Pach, and Spencer~\cite{ErdosGoldbergPachSpencer1988} extended this result to graphs with $e(G)=p\binom{n}{2}$ edges. They showed that
\begin{equation*}
\max_{S\subseteq V(G)}\left|e_G(S)-p\binom{|S|}{2}\right|\ge c\sqrt{\min\{p,1-p\}}\,n^{3/2}
\end{equation*}
whenever $2/(n-1)<p<1-2/(n-1)$. They also studied discrepancy for subsets containing half of the vertices. Erd\H{o}s, Faudree, Rousseau, and Schelp~\cite{ErdosFaudreeRousseauSchelp1994} studied local edge-density conditions that force triangles. Questions about uniform edge distribution also arise in the study of quasi-random graphs by Chung, Graham, and Wilson~\cite{ChungGrahamWilson1989} and in the work of Mubayi and R\"odl~\cite{MubayiRodl2004} on uniformity in hypergraphs. General accounts of discrepancy can be found in Beck and S\'os~\cite{BeckSos1995} and Chazelle~\cite{Chazelle2000}.

Bollob\'as and Scott~\cite{BollobasScott2006} studied discrepancy in graphs and hypergraphs. In particular, they proved that every graph $G$ with $e(G)=p\binom{n}{2}$ and $p(1-p)\ge 1/n$ satisfies
\begin{equation*}
\operatorname{disc}^{+}(G)\operatorname{disc}^{-}(G)\ge cp(1-p)n^3,
\end{equation*}
where $c>0$ is an absolute constant. In their later work on graph intersections, Bollob\'as and Scott~\cite{BollobasScott2011} proved that, for $16/n\le p,q\le 1-16/n$,
\begin{equation*}
\operatorname{disc}^{+}(G,H)\operatorname{disc}^{-}(G,H)\ge\frac{p^4(1-p)^4q^4(1-q)^4}{10^{20}}n^3
\end{equation*}
where $\operatorname{disc}^{+}(G,H)$ and $\operatorname{disc}^{-}(G,H)$ denote the largest positive and negative deviations of the overlap from its mean, respectively. As a consequence, they confirmed a conjecture of Erd\H{o}s, Goldberg, Pach, and Spencer~\cite{ErdosGoldbergPachSpencer1988} on bipartite discrepancy. They also raised the following $L^2$ problem.
\begin{problem}[\textbf{Bollob\'as and Scott}~\cite{BollobasScott2011}]
	Let \(G\) and \(H\) be graphs with \(e(G)=p\binom{n}{2}\) and \(e(H)=q\binom{n}{2}\), and let \(\pi \in S_n\) be chosen uniformly at random. Determine
	\begin{equation*}
		\min_{G,H}\left(\mathbb{E}_\pi\left( |E(G_{\pi})\cap E(H)| - pq\binom{n}{2}\right)^2\right)^{1/2}.
	\end{equation*}
\end{problem}
 Following Bollob\'as and Scott~\cite{BollobasScott2011}, we associate $G$ and $H$ with vectors $g,h\in\mathbb R^{\binom{n}{2}}$, defined 
\begin{equation*}
g(ij)=\mathbf{1}_{\{ij\in E(G)\}}-p,\qquad h(ij)=\mathbf{1}_{\{ij\in E(H)\}}-q.
\end{equation*}
For $\pi\in S_n$, let $g_\pi$ denote the vector corresponding to the relabelled graph $G_\pi$. Then
\begin{equation}\label{eq:<g,h>}
\langle g_\pi,h\rangle=|E(G_\pi)\cap E(H)|-pq\binom{n}{2}.
\end{equation}
In particular, $\mathbb E_\pi\langle g_\pi,h\rangle=0$, where $\pi$ is chosen uniformly at random from $S_n$. Our first result resolves this problem.
\begin{theorem}\label{thm:asymptotic-intro}
Let $G$ and $H$ be graphs of order $n$, and suppose that $e(G)=p\binom{n}{2}$ and $e(H)=q\binom{n}{2}$, where $p,q\in(0,1)$ are fixed. Then
\begin{equation}\label{eq:asymptotic-intro}
\min_{\substack{G,H}}\left(\mathbb E_\pi \langle g_\pi,h\rangle^2\right)^{1/2}=\left(\sqrt{\frac{p(1-p)q(1-q)}{2}\cdot\frac{1-\sqrt R}{1+\sqrt R}}\right)n+o(n),
\end{equation}
where $R=\max\{p,1-p,q,1-q\}$.
\end{theorem}
The leading coefficient is attained asymptotically by taking one graph almost regular and the other quasi-complete or quasi-star.

We next consider intersections of uniform hypergraphs. Bollob\'as and Scott~\cite{BollobasScott2015} studied intersections of hypergraphs, and their work on random hypergraphs and tournaments~\cite{BollobasScott2015Random} considered the overlap of randomly chosen structures. Ma, Naves, and Sudakov~\cite{MaNavesSudakov2015} also obtained discrepancy results for random graphs and hypergraphs. Here we consider a particular extremal problem posed
by Bollob\'as and Scott~\cite{BollobasScott2011}.

 Let $n\ge 2$ be even, and let $\mathcal X$ be the family of all $(n/2)$-uniform hypergraphs on $n$ vertices with exactly $\frac12\binom{n}{n/2}$ edges.
For $\mathcal{G},\mathcal{H}\in\mathcal X$, define
\begin{equation*}
\operatorname{disc}^{+}(\mathcal{G},\mathcal{H})=\max_{\pi\in S_n}|E(\mathcal{G}_\pi)\cap E(\mathcal{H})|-\frac14\binom{n}{n/2},
\end{equation*}
and
\begin{equation*}
\operatorname{disc}^{-}(\mathcal{G},\mathcal{H})=\frac14\binom{n}{n/2}-\min_{\pi\in S_n}|E(\mathcal{G}_\pi)\cap E(\mathcal{H})|.
\end{equation*}
\begin{problem}[\textbf{Bollob\'as and Scott}~\cite{BollobasScott2011}]
Consider the family \(\mathcal{X}\) of \((n/2)\)-uniform hypergraphs on \(n\) vertices with \(\frac{1}{2}\binom{n}{n/2}\) edges. Determine
\begin{equation*}
	\min_{\mathcal{G},\mathcal{H}\in\mathcal X}\operatorname{disc}^{+}(\mathcal{G},\mathcal{H})\operatorname{disc}^{-}(\mathcal{G},\mathcal{H}).
\end{equation*}
\end{problem}
Our second result resolves this problem.
\begin{theorem}\label{thm:hypergraph-intro}
Let $n\ge 2$ be even, and let $\mathcal X$ be the family of all $(n/2)$-uniform hypergraphs on $n$ vertices with exactly $\frac12\binom{n}{n/2}$ edges. If \(4 \mid \binom{n}{n/2}\), then
\begin{equation*}\label{eq:discrepancy-values-intro2}
	\min_{\mathcal{G},\mathcal{H}\in\mathcal X}\operatorname{disc}^{+}(\mathcal{G},\mathcal{H})=
	\min_{\mathcal{G},\mathcal{H}\in\mathcal X}\operatorname{disc}^{-}(\mathcal{G},\mathcal{H})=0,
\end{equation*}
Otherwise, 
\begin{equation*}\label{eq:discrepancy-values-intro}
\min_{\mathcal{G},\mathcal{H}\in\mathcal X}\operatorname{disc}^{+}(\mathcal{G},\mathcal{H})=
	\min_{\mathcal{G},\mathcal{H}\in\mathcal X}\operatorname{disc}^{-}(\mathcal{G},\mathcal{H})=\frac12.
\end{equation*}
Consequently,
\begin{equation*}
\min_{\mathcal G,\mathcal H\in\mathcal X}
\operatorname{disc}^{+}(\mathcal G,\mathcal H)
\operatorname{disc}^{-}(\mathcal G,\mathcal H)
=\begin{cases}
	0, &  4\mid\binom{n}{n/2},\\
	\frac14, & \text{otherwise}.
\end{cases}
\end{equation*}
\end{theorem}
The paper is organized as follows. Section~\ref{sec:Auxiliary lemmas} presents the preliminary results needed for the graph problem. Section~\ref{sec:The minimum $L^2$ norm} proves Theorem~\ref{thm:asymptotic-intro}. Section~\ref{sec:Discrepancy of ${n}/{2}$-uniform hypergraphs} proves Theorem~\ref{thm:hypergraph-intro}.
\section{Auxiliary lemmas}\label{sec:Auxiliary lemmas}
In this section, we establish some lemmas on graphs and hypergraphs that will be used later.

For a graph $G$ of order $n\geq 4$ with $e(G)=p\binom{n}{2}$, define
\begin{equation}\label{eq:AG,Tp}
A(G) = \frac{1}{n-2} \sum_{v \in V(G)} \left( d_G(v) - \frac{2e(G)}{n} \right)^2, \qquad T_{p} = p(1-p)\binom{n}{2}.
\end{equation}
For a graph $H$ on $n\geq 4$ with $e(H)=q\binom{n}{2}$, define $A(H)$ and $T_{q}$ analogously.
\begin{lemma}\label{lem:second-moment}
Let $G$ and $H$ be graphs of order $n\ge4$, and suppose that $e(G)=p\binom{n}{2}$ and $e(H)=q\binom{n}{2}$. Then
	\begin{equation*}\label{eq:second-moment}
\mathbb E_\pi\langle g_\pi,h\rangle^2=\frac{A(G)A(H)}{n-1}+\frac{(T_p-A(G))(T_q-A(H))}{\binom{n}{2}-n}.
	\end{equation*}
\end{lemma}
\begin{proof}[\bf Proof]
For graphs $G$ and $H$, let
\begin{equation}\label{eq:pGpH}
P_G:=\sum_{v\in V(G)}\binom{d_G(v)}{2},\qquad P_H:=\sum_{v\in V(H)}\binom{d_H(v)}{2}.
\end{equation}
Thus $P_G$ is the number of unordered pairs of adjacent edges in $G$, while $\binom{e(G)}{2}-P_G$ is the number of unordered pairs of disjoint edges. The complete graph $K_n$ contains $n\binom{n-1}{2}$ unordered pairs of adjacent edges and $3\binom{n}{4}$ unordered pairs of disjoint edges. Hence
\begin{align*}
	\mathbb E_\pi |E(G_\pi)\cap E(H)|^2&=\mathbb{E}_{\pi}|E(G_{\pi}) \cap E(H)| + 2\mathbb{E}_{\pi}\binom{|E(G_{\pi}) \cap E(H)|}{2}\\&=pq\binom{n}{2}+\frac{2P_GP_H}{n\binom{n-1}{2}}+\frac{2\left(\binom{e(G)}{2}-P_G\right)\left(\binom{e(H)}{2}-P_H\right)}{3\binom{n}{4}}.
\end{align*}
By \eqref{eq:<g,h>} and $\mathbb E_\pi\langle g_\pi,h\rangle=0$, we have
\begin{equation*}
\mathbb E_\pi |E(G_\pi)\cap E(H)|=pq\binom{n}{2}.
\end{equation*}
Then
\begin{align}\label{eq:E<g,h>^2}
\mathbb E_\pi\langle g_\pi,h\rangle^2
&=\mathbb E_\pi\left(|E(G_\pi)\cap E(H)|-pq\binom{n}{2}\right)^2\notag
\\&=\mathbb E_\pi |E(G_\pi)\cap E(H)|^2-2pq\binom{n}{2}\,\mathbb E_\pi |E(G_\pi)\cap E(H)|+p^2q^2\binom{n}{2}^{2}\notag
\\&=\mathbb E_\pi |E(G_\pi)\cap E(H)|^2-p^2q^2\binom{n}{2}^{2}\notag
\\&=pq\binom{n}{2}+\frac{2P_GP_H}{n\binom{n-1}{2}}+\frac{2\left(\binom{e(G)}{2}-P_G\right)\left(\binom{e(H)}{2}-P_H\right)}{3\binom{n}{4}}-p^2q^2\binom{n}{2}^{\,2}.
\end{align}
Since
\begin{equation}\label{eq:PG}
P_G=\frac12\left(\sum_{v\in V(G)}d_G(v)^2-2e(G)\right)=\frac12\left((n-2)A(G)+\frac{4e(G)^2}{n}-2e(G)\right),
\end{equation}
with the analogous identity holding for $P_H$. Combining \eqref{eq:AG,Tp}, \eqref{eq:E<g,h>^2} and \eqref{eq:PG}, we have
\begin{align*}
\mathbb E_\pi\langle g_\pi,h\rangle^2
&=\frac{A(G)A(H)}{n-1}+\frac{A(G)A(H)-A(G)T_q-T_pA(H)+T_pT_q}{\binom{n}{2}-n}
\\&=\frac{A(G)A(H)}{n-1}+\frac{(T_p-A(G))(T_q-A(H))}{\binom{n}{2}-n}.
\end{align*}
This completes the proof.
\end{proof}
We next give a simple lemma on vertex degrees, which will be used to bound $A(G)$ and to prove the hypergraph result.
\begin{lemma}\label{lem:balanced-hypergraph}
	For integers $1\le r\le n$ and $0\le m\le\binom{n}{r}$, there exists a simple $r$-uniform
	hypergraph $H$ on $[n]$ with $m$ edges such that $|d_H(u)-d_H(v)|\le1$ for all $u,v\in[n]$.
\end{lemma}
\begin{proof}[\bf Proof]
	Choose a simple $r$-uniform hypergraph $H$ on $[n]$ with $m$ edges minimizing $\sum_{w\in[n]} d_H(w)^2$. Suppose that $d_H(u)\ge d_H(v)+2$ for some $u,v\in[n]$. Then
	\begin{equation*}
|\{E\in E(H):u\in E,\ v\notin E\}|-|\{E\in E(H):v\in E,\ u\notin E\}|=d_H(u)-d_H(v)>0.
	\end{equation*}
Since replacing $u$ by $v$ in distinct edges gives distinct sets, there is an edge $E\in E(H)$ with $u\in E$, $v\notin E$ and $(E\setminus\{u\})\cup\{v\}\notin E(H)$. Replacing $E$ by $(E\setminus\{u\})\cup\{v\}$, we obtain
\begin{equation*}
(d_H(u)-1)^2+(d_H(v)+1)^2-d_H(u)^2-d_H(v)^2=-2(d_H(u)-d_H(v))+2<0,
\end{equation*}
which contradicts the choice of $H$. Hence, $|d_H(u)-d_H(v)|\le1$ for all $u,v\in[n]$.
\end{proof}
\begin{remark}\label{remark:dv-du<=1}
For a graph $G$ with a fixed number of vertices and edges that minimizes $A(G)$, we have $|d_G(u)-d_G(v)|\le 1$ for all $u,v\in V(G)$.
\end{remark}
The remark allows us to determine the minimum of $A(G)$ exactly.
\begin{lemma}\label{lem:min-variance}
	Let $G$ be a graph of order $n\ge4$ for which $A(G)$ is minimum among all graphs of order $n$ with $e(G)$ edges. Write $2e(G)=an+\sigma$ with $0\le \sigma<n$. Then
	\begin{equation*}
	A(G)=\frac{\sigma(n-\sigma)}{n(n-2)}.
	\end{equation*}
Moreover, $G$ has $\sigma$ vertices of degree $a+1$ and $n-\sigma$ vertices of degree $a$.
\end{lemma}
\begin{proof}[\bf Proof]
	Since
	\begin{align}\label{eq:A(G)2}
A(G)&=\frac{1}{n-2}\sum_{v\in V(G)}\left(d_G(v)-\frac{2e(G)}{n}\right)^2\notag
\\&=\frac{1}{n-2}\left(\sum_{v\in V(G)}d_G(v)^2-\frac{4e(G)}{n}\sum_{v\in V(G)}d_G(v)+\frac{4e(G)^2}{n}\right)\notag
\\&=\frac{1}{n-2}\left(\sum_{v\in V(G)}d_G(v)^2-\frac{4e(G)^2}{n}\right),
	\end{align}
	and $e(G)$ is fixed, minimizing $A(G)$ is equivalent to minimizing $\sum_{v\in V(G)}d_G(v)^2$. By Remark~\ref{remark:dv-du<=1}, we claim that any two degrees in $G$ differ by at most one. Since
	\begin{equation*}
		\sum_{v\in V(G)}d_G(v)=2e(G)=an+\sigma
	\end{equation*}
	and $0\le\sigma<n$, exactly $\sigma$ vertices have degree $a+1$, while the remaining $n-\sigma$ vertices have degree $a$. The average degree is ${2e(G)}/{n}=a+{\sigma}/{n}$. Using \eqref{eq:A(G)2}, we have
	\begin{equation*}
	A(G)=
	\frac{1}{n-2}
	\left[
	\sigma\left(1-\frac{\sigma}{n}\right)^2
	+(n-\sigma)\left(\frac{\sigma}{n}\right)^2
	\right]=
	\frac{\sigma(n-\sigma)}{n(n-2)}.
	\end{equation*}
This completes the proof.
\end{proof}
We next determine the maximum of $A(G)$ for fixed $n$ and $e(G)$. Write
\begin{equation*}
\binom{t}{2}\le e(G)<\binom{t+1}{2},\qquad e(G)=\binom{t}{2}+\nu,
\end{equation*}
where $0\le \nu<t$. A quasi-complete graph with $e(G)$ edges consists of $K_t$, one additional vertex adjacent to $\nu$ vertices of $K_t$ if $\nu>0$, and isolated vertices. 
\begin{equation}\label{eq:sumdv^2}
	\sum_{v\in V(G)}d_G(v)^2=\nu t^2+(t-\nu)(t-1)^2+\nu^2.
\end{equation}
A quasi-star graph with $e(G)$ edges is the complement of a quasi-complete graph with $\binom{n}{2}-e(G)$ edges. Since $d_{\overline G}(v)=n-1-d_G(v)$,
\begin{equation*}
	\sum_{v\in V(\overline G)}d_{\overline G}(v)^2=\sum_{v\in V(G)}(n-1-d_G(v))^2=n(n-1)^2-4(n-1)e(G)
	+\sum_{v\in V(G)}d_G(v)^2.
\end{equation*}
Thus the squared-degree sum of a quasi-star follows from \eqref{eq:sumdv^2} by replacing $e(G)$ with $\binom{n}{2}-e(G)$ and taking the complement. Since
\begin{equation*}
	A(G)=\frac{1}{n-2}\left(\sum_{v\in V(G)}d_G(v)^2-\frac{4e(G)^2}{n}\right),
\end{equation*}
and
\begin{equation*}
\sum_{v\in V(G)}\binom{d_G(v)}{2}=\frac12\left(\sum_{v\in V(G)}d_G(v)^2-2e(G)\right),
\end{equation*}
maximizing $A(G)$ for fixed $n$ and $e(G)$ is equivalent to maximizing the number of unordered pairs of adjacent edges. The theorem of Ahlswede and Katona~\cite{AhlswedeKatona1978} therefore gives the following result.
\begin{theorem}[\textbf{Ahlswede and Katona}~\cite{AhlswedeKatona1978}]\label{thm:max-variance}
	For graphs with fixed numbers of vertices and edges, $A(G)$ is maximized by a quasi-complete graph or a quasi-star graph.
\end{theorem}
\section{The minimum $L^2$ norm}\label{sec:The minimum $L^2$ norm}
In this section, we use the results of Section~\ref{sec:Auxiliary lemmas} to prove Theorem~\ref{thm:asymptotic-intro}.
\begin{lemma}\label{lem:endpoint-asymptotics}
Fix $p\in(0,1)$, and let $G$ be range over graphs on $n$ vertices with $e(G)=p\binom{n}{2}$. Then
\begin{equation*}
\min_{G} A(G) = O(1)
\end{equation*}
and
\begin{equation}\label{eq:Tp-maxAG}
T_p - \max_{G} A(G) = \frac{1}{2}p(1-p)\frac{1-\sqrt{\max\{p,1-p\}}}{1+\sqrt{\max\{p,1-p\}}}n^2 + O_p(n).
\end{equation}
\end{lemma}
\begin{proof}[\bf Proof]
By Lemma~\ref{lem:min-variance}, we have 
\begin{equation*}
\min_G A(G)=\frac{\sigma(n-\sigma)}{n(n-2)}\le\frac{n}{4(n-2)}=O(1),
\end{equation*}
where $2e(G)=an+\sigma$ and $0\le\sigma<n$. By Theorem~\ref{thm:max-variance}, the maximum of $A(G)$ is attained by a quasi-complete graph or a quasi-star graph. We estimate these two cases separately.

\item \textbf{Case 1.}  Let $G$ be a quasi-complete graph with $p\binom{n}{2}$ edges.
\begin{equation*}
\binom{t}{2}\le p\binom{n}{2}<\binom{t+1}{2},\qquad p\binom{n}{2}=\binom{t}{2}+\nu,
\end{equation*}
where $0\le\nu<t$. It follows that $t=\sqrt p\,n+O_p(1)$ and $\nu=O_p(n)$. By \eqref{eq:sumdv^2}, we have
\begin{equation*}
\sum_{v\in V(G)}d_G(v)^2=t(t-1)^2+\nu(2t-1)+\nu^2=t^3+O_p(n^2)=p^{3/2}n^3+O_p(n^2).
\end{equation*}
Therefore, by \eqref{eq:A(G)2}, we have
\begin{equation*}
A(G)=\frac{1}{n-2}\left(p^{3/2}n^3-\frac{4p^2\binom{n}{2}^2}{n}+O_p(n^2)\right)=\bigl(p^{3/2}-p^2\bigr)n^2+O_p(n).
\end{equation*}
Since
\begin{equation*}
T_p=p(1-p)\binom{n}{2}=\frac12p(1-p)n^2+O_p(n),
\end{equation*}
we obtain
\begin{equation}\label{eq:T_p-A(G)}
T_p-A(G)=\left(\frac12p(1-p)-p^{3/2}+p^2\right)n^2+O_p(n)=\frac12p(1-\sqrt p)^2n^2+O_p(n).
\end{equation}
\item \textbf{Case 2.} Let $G$ be a quasi-star graph with $p\binom{n}{2}$ edges. 

Its complement is a quasi-complete graph with $(1-p)\binom{n}{2}$ edges. Since
\begin{equation*}
d_{\overline G}(v)-\frac{2e(\overline G)}{n}=-\left(d_G(v)-\frac{2e(G)}{n}\right),
\end{equation*}
we have $A(\overline G)=A(G)$. Together with $T_p=T_{1-p}$, applying \eqref{eq:T_p-A(G)} with $p$ replaced by $1-p$ gives
\begin{equation*}
T_p-A(G)=\frac12(1-p)(1-\sqrt{1-p})^2n^2+O_p(n).
\end{equation*}
Finally, we consider
\begin{equation*}
p(1-\sqrt p)^2=p(1-p)\frac{1-\sqrt p}{1+\sqrt p},\qquad
(1-p)(1-\sqrt{1-p})^2=p(1-p)\frac{1-\sqrt{1-p}}{1+\sqrt{1-p}}.
\end{equation*}
Since ${(1-\sqrt z)}/{(1+\sqrt z)}$ is decreasing for $0<z<1$, the smaller of the two leading terms is obtained by taking $z=\max\{p,1-p\}$. By Theorem~\ref{thm:max-variance}, we have
\begin{equation*}
T_p-\max_G A(G)=\frac12p(1-p)\frac{1-\sqrt{\max\{p,1-p\}}}{1+\sqrt{\max\{p,1-p\}}}\,n^2+O_p(n).
\end{equation*}
This completes the proof.
\end{proof}
We now combine the preceding results to prove Theorem~\ref{thm:asymptotic-intro}.
\begin{proof}[\bf Proof of Theorem~\ref{thm:asymptotic-intro}]
By Lemma~\ref{lem:endpoint-asymptotics}, we have 
\begin{equation*}
\min_G A(G)=O(1),\qquad \min_H A(H)=O(1),
\end{equation*}
and
\begin{equation*}
T_p=\frac12p(1-p)n^2+O_p(n),\qquad T_q=\frac12q(1-q)n^2+O_q(n).
\end{equation*}
By Lemma~\ref{lem:second-moment}, we have
\begin{equation*}
\mathbb E_\pi\langle g_\pi,h\rangle^2=\frac{A(G)A(H)}{n-1}+\frac{(T_p-A(G))(T_q-A(H))}{\binom{n}{2}-n}.
\end{equation*}
Since this expression is affine in each of $A(G)$ and $A(H)$ when the other is fixed, its minimum is attained when each of $A(G)$ and $A(H)$ is either as small or as large as possible. 

\item \textbf{Case 1.} If both $A(G)$ and $A(H)$ attain their minimum values, then
\begin{align*}
\mathbb E_\pi\langle g_\pi,h\rangle^2
&=\frac{O(1)}{n-1}+\frac{\left(p(1-p)n^2/2+O_p(n)\right)\left(q(1-q)n^2/2+O_q(n)\right)}
{n^2/2+O(n)}
\\&=\frac12p(1-p)q(1-q)n^2+O_{p,q}(n).
\end{align*}
\item \textbf{Case 2.} If $A(G)$ is minimum and $A(H)$ is maximum, then the first term is $O_{p,q}(n)$, while Lemma~\ref{lem:endpoint-asymptotics} gives
\begin{equation*}
\mathbb E_\pi\langle g_\pi,h\rangle^2=\frac12p(1-p)q(1-q)\frac{1-\sqrt{\max\{q,1-q\}}}
{1+\sqrt{\max\{q,1-q\}}}\,n^2+O_{p,q}(n).
\end{equation*}
\item \textbf{Case 3.} If $A(G)$ is maximum and $A(H)$ is minimum, then
\begin{equation*}
\mathbb E_\pi\langle g_\pi,h\rangle^2=\frac12p(1-p)q(1-q)\frac{1-\sqrt{\max\{p,1-p\}}}
{1+\sqrt{\max\{p,1-p\}}}\,n^2+O_{p,q}(n).
\end{equation*}
\item \textbf{Case 4.} If both $A(G)$ and $A(H)$ are maximum, then $A(G)=\Theta_p(n^2)$ and $A(H)=\Theta_q(n^2)$. Hence $A(G)A(H)/(n-1)=\Theta_{p,q}(n^3)$, so this case cannot give the minimum. Therefore
\begin{equation*}
\min_{G,H}\mathbb E_\pi\langle g_\pi,h\rangle^2=\frac12p(1-p)q(1-q)\frac{1-\sqrt R}{1+\sqrt R}\,n^2+O_{p,q}(n),
\end{equation*}
where $R=\max\{p,1-p,q,1-q\}$. 

Since the leading coefficient is positive for fixed $p,q\in(0,1)$, we obtain
\begin{align*}
\min_{G,H}\left(\mathbb E_\pi\langle g_\pi,h\rangle^2\right)^{1/2}
&=\sqrt{\frac{p(1-p)q(1-q)}{2}\frac{1-\sqrt R}{1+\sqrt R}}\,n+O_{p,q}(1)
\\&=\sqrt{\frac{p(1-p)q(1-q)}{2}\frac{1-\sqrt R}{1+\sqrt R}}\,n+o(n),
\end{align*}
which proves the theorem.
\end{proof}
\section{Discrepancy of ${n}/{2}$-uniform hypergraphs}\label{sec:Discrepancy of ${n}/{2}$-uniform hypergraphs}
In this section, we prove Theorem~\ref{thm:hypergraph-intro} by giving the required bounds and constructions.
\begin{proof}[\bf Proof of Theorem~\ref{thm:hypergraph-intro}]
For any $\mathcal{G},\mathcal{H}\in\mathcal X$, each edge of $\mathcal H$ belongs to $\mathcal G_\pi$ with probability $1/2$, where $\pi\in S_n$ is chosen uniformly at random. Thus
\begin{equation}\label{eq:hypergraph-mean}
	\mathbb E_\pi |E(\mathcal{G}_\pi)\cap E(\mathcal{H})|=\frac12 e(\mathcal{H})=\frac14\binom{n}{n/2}.
\end{equation}
Since $\binom{n}{n/2}=2\binom{n-1}{n/2-1}$ is even, we consider the following two cases.
\item \textbf{Case 1.} $4\mid\binom{n}{n/2}$.

Let $\mathcal G^*$ be the full star consisting of all $(n/2)$-subsets containing a fixed vertex. Since 
\begin{equation*}
e(\mathcal G^*)=\binom{n-1}{n/2-1}=\frac12\binom{n}{n/2},
\end{equation*}
we have $\mathcal G^*\in\mathcal X$. By Lemma~\ref{lem:balanced-hypergraph}, choose $\mathcal H\in\mathcal X$ such that $|d_{\mathcal H}(u)-d_{\mathcal H}(v)|\le1$ for all $u,v\in[n]$. Then
\begin{equation*}
d_{\mathcal{H}}=\frac1n\sum_{v=1}^n d_{\mathcal H}(v)=\frac{(n/2)e(\mathcal H)}n=\frac14\binom{n}{n/2}.
\end{equation*}
Since the average degree is an integer, we obtain
\begin{equation*}
d_{\mathcal H}(v)=\frac14\binom{n}{n/2}
\end{equation*}
for every $v\in[n]$. For each $\pi\in S_n$, $\mathcal G_\pi^*$ is again a full star. Let $\mathcal H_\pi^*$ be the subhypergraph of $\mathcal H$ consisting of all edges containing the center of $\mathcal G_\pi^*$. Then
\begin{equation*}
|E(\mathcal G_\pi^*)\cap E(\mathcal H)|=e(\mathcal H_\pi^*)=\frac14\binom{n}{n/2}
\end{equation*}
for every $\pi\in S_n$. Therefore $\operatorname{disc}^{+}(\mathcal G^*,\mathcal H)=\operatorname{disc}^{-}(\mathcal G^*,\mathcal H)=0$.

\item \textbf{Case 2.} $4\nmid\binom{n}{n/2}$.

Since $\binom{n}{n/2}$ is even, we have $\binom{n}{n/2}\equiv2\pmod4$. The mean in \eqref{eq:hypergraph-mean} is a half-integer, whereas $|E(\mathcal G_\pi)\cap E(\mathcal H)|$ is an integer for every $\pi\in S_n$. Thus
\begin{equation*}
\min_{\pi\in S_n}|E(\mathcal G_\pi)\cap E(\mathcal H)|\le\frac14\binom{n}{n/2}-\frac12,\qquad
\max_{\pi\in S_n}
|E(\mathcal G_\pi)\cap E(\mathcal H)|
\ge\frac14\binom{n}{n/2}+\frac12.
\end{equation*}
Therefore, for every $\mathcal G,\mathcal H\in\mathcal X$,
\begin{equation*}
\operatorname{disc}^{+}(\mathcal G,\mathcal H)\ge\frac12,\qquad\operatorname{disc}^{-}(\mathcal G,\mathcal H)\ge\frac12.
\end{equation*}
To attain these bounds, take $\mathcal G^*$ as in Case~1 and choose $\mathcal H\in\mathcal X$ such that $|d_{\mathcal H}(u)-d_{\mathcal H}(v)|\le1$ for all $u,v\in[n]$. Since the average degree is $\frac14\binom{n}{n/2}$, which is a half-integer,
\begin{equation*}
d_{\mathcal H}(v)\in\left\{\frac14\binom{n}{n/2}-\frac12,\,\frac14\binom{n}{n/2}+\frac12\right\}
\end{equation*}
for every $v\in[n]$. As $\pi$ ranges over $S_n$, the center of $\mathcal G_\pi^*$ ranges over all vertices of $[n]$. Hence
\begin{equation*}
\min_{\pi\in S_n}|E(\mathcal G_\pi^*)\cap E(\mathcal H)|=\min_{\pi\in S_n}e(\mathcal H_\pi^*)=\frac14\binom{n}{n/2}-\frac12,
\end{equation*}
and
\begin{equation*}
\max_{\pi\in S_n}|E(\mathcal G_\pi^*)\cap E(\mathcal H)|=\max_{\pi\in S_n}e(\mathcal H_\pi^*)=\frac14\binom{n}{n/2}+\frac12.
\end{equation*}
Consequently, $\operatorname{disc}^{+}(\mathcal G^*,\mathcal H)=\operatorname{disc}^{-}(\mathcal G^*,\mathcal H)=1/2$. This completes the proof.
\end{proof}
\section*{Declaration on the use of AI}
The author used ChatGPT 5.6 Pro to assist in discussing proof strategies, checking proofs, and improving exposition.

\end{document}